\documentclass[A4paper,11pt]{article}
\usepackage{latexsym}
\usepackage{mathrsfs}
\usepackage{amssymb}
\usepackage{amscd}
\usepackage[dvips]{graphicx}                  

\newtheorem{theorem}{Theorem}

\newtheorem{lemma}{Lemma}

\newenvironment{proof}[1][Proof]{\textbf{#1.} }{\ \rule{0.5em}{0.5em}}
\date{}
\long\def\symbolfootnote[#1]#2{\begingroup%
	\def\thefootnote{$\;$}\footnote[#1]{$^*$#2}\endgroup}
\begin{document}
	
	\title{Chains in the Rudin-Frol\'ik order for regulars}
	\author{Joanna Jureczko}
\maketitle

\symbolfootnote[2]{Mathematics Subject Classification: Primary 03E10, 03E20, 03E30..

\hspace{0.2cm}
Keywords: \textsl{Ultrafilters, regular cardinal, Rudin-Frol\'ik order, independent family.}}

\begin{abstract}
	The aim of this paper is to construct chains of length $(2^\kappa)^+$ in the Rudin-Frol\'ik order of $\beta \kappa$ for $\kappa$ regular.
\end{abstract}

\section{Introduction} 

The investigations of the Rudin-Frol\'ik order, which has been defined by Z. Frol\'ik in \cite{ZF}, are an important topic, but still little known. Initially, this order was used by Frol\'ik   to prove that $\beta\omega\setminus \omega$ is not homogeneous.  M.E. Rudin, who nearly defined this ordering  in \cite{MR1},  as the first observed that the relation between filters she used is really  ordering.
Further this order was investigated by D. Booth in \cite{DB} who showed that this relation is a partial ordering of the equivalence classes, that is a tree, and that it is not well-founded.

Later, in \cite{MR} the author defined and studied the partial orders on the type of points in $\beta\omega$ and in $\beta\omega \setminus \omega$. These definitions were used later in \cite{BB} and \cite{BE3}.

The papers which are included in the scope of our considerations are E. Butkovi\v cov\'a's results who worked on this topic between 1981 and 1990 publishing a number of papers concerning ultrafilters  in the Rudin-Frol\'ik order in $\beta \omega\setminus \omega$.  
Let us briefly review her achievements in this topic.
In \cite{BB} with L. Bukovsk\'y and in \cite{BE3}, she constructed an ultrafilter on $\omega$ with the countable set of its predecessors. In \cite{BE1} she  constructed ultrafilters without immediate predecessors. In \cite{BE3}, Butkovi\v cov\'a showed that there exists in the Rudin-Frol\'ik order an unbounded chain  orded-isomorphic in $\omega_1$.  In \cite{BE4}, she proved that  there is a set of $2^{2^{\aleph_0}}$ ultrafilters incomparable in the Rudin-Frol\'ik order which is bounded from below and no its subset of cardinality more than one has an infimum. In \cite{BE5}, Butkovi\v cov\'a proved that for every cardinal between $\omega$ and $\mathfrak{c}$ there is a strictly decreasing chain without a lower bound.
In most of these papers there is used method presented in \cite{KK}.

In 1976, A. Kanamori published a paper \cite{AK} in which, among others showed that the Rudin-Frol\'ik tree cannot be very high if one considers it over a measurable cardinal. Moreover, in the same paper he left a number of open problems about the Rudin-Frol\'ik order. 
Recently,  M. Gitik in \cite{MG} answered some of them but using metamathematical methods. The solution of some of the problems from \cite{AK} presented in combinatorial methods are in preparation, (\cite{JJ_kanamori}).

However, the Rudin-Frol\'ik order was investigated mainly for $\beta\omega$, significant results may be obtained when considering this order for the space $ \beta \kappa $, where $ \kappa $ is any cardinal. 
Since some of the results are based on the construction of the sequences of filters  by transfiite induction, a different technique is needed in the case of $\beta \kappa$ when $\kappa > \omega$.

The method proposed in \cite{BK} by Baker and Kunen  comes in handy.
In mentioned paper the authors presented  very usefull method which can be recognized as a generalization of method  presented in \cite{KK}. 
It is worth empahsizing that both methods, (from \cite{KK} and \cite{BK}), provide usefull "technology"  for keeping the transfinite construction for an ultrafilter not finished before $\mathfrak{c}$ steps, (see \cite{KK}), and $2^\kappa, $  for $\kappa$ being infinite cardinal, (see \cite{BK}), but the second method has some limitations, among others $\kappa$ must be regular.
Due to the lack of adequate useful method for a singular cardinal $\kappa$, the similar results but for singulars are still left as open questions. So far, we have not found an answer whether the assumptions can be omitted, which would probably also involve  changes in the methods used in our considerations. Therefore, based on the results from \cite{BK}, we restrict our results to this particular case.

The results contained in this paper are a continuation of research on the properties of subsets of $\beta\kappa$ (\cite{JJ_order, JJ_order2}), but due to the methods used here and taken from \cite{BK}, they are also limited to the case where $\kappa$  is a regular cardinal.

The starting point of investigations presented here are papers \cite{BE3} and \cite{JJ_order}. In \cite{BE3} there are constructed chains of ultrafilters in the Rudin-Frol\'ik order of $\beta \omega$ order-isomorphic to $(2^{\aleph_0})^+$, while in \cite{JJ_order} there are presented constructions of ultrafilters in $\beta \kappa$  in the Rudin Frol\'ik order, for regular $\kappa$.
As was mentioned in \cite{BE3}, there are at most $2^{\aleph_0}$ predecessors of the type of an ultrafilter, and cardinality of each branch in the  Rudin-Frol\'ik order is at least $2^{\aleph_0}$. Therefore, it would expect that long chains of ultrafilters in $\beta \kappa$ for regular $\kappa$  can have cardinality  $2^\kappa$  or $(2^{\kappa})^+$. In this paper, there is constructed a chain of ultrafilters having  cardinality $(2^\kappa)^+$. 

Since the main method used in the results presented in this paper is based on the method proposed in \cite{BK}, we only show  case $(2^\kappa)^+$, for regular $\kappa$. We still do not know what is for singulars and we leave it as an open problem.
Note that the considerations have further consequences which are given in \cite{JJ_order, JJ_order2, JJ_order3}.
The results included in this paper were archived for the first time as the report of Wroclaw University of Science ant Technology on 2.12.2022, \cite{JJ_order4}.

The paper is organized as follows: in Section 2, there are presented definitions and previous facts needed for the results presented in further parts of this paper. In Section 3, there are proved auxiliary results, the main result and the open problem.

We have tried to present all the necessary definitions, assuming tacitly that the reader has a basic knowledge of ultrafilters and the Rudin-Frol\'ik order.
However, for definitions and facts not quoted here, we refer the reader to e.g. \cite{CN, TJ}.

\section{Definitions and previous results}

   \textbf{2.1.}
   In the whole paper, we assume that $\kappa$ is an infinite cardinal. Then $\beta\kappa$ means the \v Cech-Stone compactification, where $\kappa$ has the discrete topology. Hence, $\beta\kappa$ is the space of ultrafilters on $\kappa$ and $\beta \kappa \setminus \kappa$ is the space of nonprincipal ultrafilters on $\kappa$. 
   \\\\
   \textbf{2.2.}
   A set $\{\mathcal{F}_\alpha \colon \alpha < \kappa\}$ of filters on $\kappa$ is \textit{$\kappa$-discrete} iff there is a partition $\{A_\alpha \colon \alpha < \kappa\}$ of $\kappa$ such that $A_\alpha \in \mathcal{F}_\alpha$ for each $\alpha < \kappa$. 
   \\\\
  \textbf{2.3.}  Let $\mathcal{F}, \mathcal{G}\in \beta \kappa\setminus \kappa$. We define the \textit{Rudin-Frol\'ik order} as follows
  $$\mathcal{F} \leqslant_{RF} \mathcal{G} \textrm{ iff } \mathcal{G} = \Sigma(X, \mathcal{F})$$
  for some $\kappa$-discrete set $X = \{\mathcal{F}_\alpha \colon \alpha < \kappa\} \subseteq \beta\kappa,$
  where $$\Sigma(X, \mathcal{F}) = \{A \subseteq \kappa \colon \{\alpha < \kappa \colon A \in \mathcal{F}_\alpha\}\in \mathcal{F}\}.$$
  We define $$\mathcal{F} =_{RF} \mathcal{G} \textrm{ iff } \mathcal{F} \leqslant_{RF} \mathcal{G} \textrm{ and } \mathcal{G} \leqslant_{RF} \mathcal{F}$$  $$\mathcal{F} <_{RF} \mathcal{G} \textrm{ iff } \mathcal{F} \leqslant_{RF} \mathcal{G} \textrm{ and } \mathcal{F} \not =_{RF} \mathcal{G}.$$
   \textbf{2.4.} Let $\tau$ and $\kappa$ be infinite cardinals. A set of filters $\{\mathcal{F}_{\xi, \zeta} \colon \xi < \tau, \zeta < \kappa\}$ is called \textit{stratified} iff
   \begin{itemize}
   	\item [(1)] $\{\mathcal{F}_{\xi, \zeta} \colon \zeta < \kappa\}$ is $\kappa$-discrete for each $\xi < \tau$,
   	\item [(2)] for each $\xi < \tau, \zeta < \kappa$ and each $\gamma$ such that $\xi < \gamma < \tau$   $$|\{\mu < \kappa \colon A \in \mathcal{F}_{\gamma, \mu}\}| = \kappa$$
   	for all $A \in \mathcal{F}_{\xi, \zeta}$. 
   \end{itemize}
   \textbf{2.5.} Let $Y =\{\mathcal{F}_{\xi, \zeta} \colon \xi < \tau, \zeta < \kappa\}$ be a stratified set of filters and let $W$ be a subset of $Y$. We define
   \begin{itemize}
   	\item [(1)] $W(0) = W$,
   	\item [(2)] $W(\gamma) = \bigcup_{\beta < \gamma} W(\beta)$ for limit $\gamma$,
   	\item [(3)] $W(\gamma+1) = W(\gamma) \cup \{\mathcal{F}_{\xi, \zeta} \colon \exists_{\eta> \gamma}\  \exists_{{A \in \mathcal{F}_{\xi, \zeta}}} \{\mathcal{F}_{\eta, \nu}\colon A \in \mathcal{F}_{\eta, \nu}\}\subseteq W(\gamma)\}$,
   	\item [(4)] $\tilde{W} = \bigcup_{\gamma < \kappa^+} W(\gamma)$.
   \end{itemize}
   
   Intuitively, the above construction is used to select only certain filters from $Y$ with the desired property and then add (inductively) to the set $ W $ only those filters outside $W$ which satisfy the condition (2) in the definition 2.4. This construction   will be used in the proof of Lemma 1, where the formulation of which would not be possible taking the entire set of ultrafiters into account.
   \\\\
\textbf{2.6.} 	Let us accept the following notation: 
\begin{itemize}
	\item $\mathcal{FR}(\kappa) = \{A \subset \kappa \colon |\kappa\setminus A|< \kappa\},$
	\item $[A, B, C,...]$ means a filter generated by $A, B, C, ...$.
\end{itemize}
   \textbf{2.7.} A function $\hat\varphi \colon [\kappa^+]^{<\omega} \to [\kappa]^{< \omega}$ is \textit{$\kappa$-shrinking} iff 
   \begin{itemize}
   	\item [(1)] $p\subseteq q$ implies $\hat{\varphi}(p) \subseteq \hat{\varphi}(q)$, for any $p, q \in [\kappa^+]^{<\omega}$,
   	\item [(2)] $\hat{\varphi}(0) = 0$.
   \end{itemize}
   
   A \textit{step-family} (over $\kappa$, with respect to $\hat{\varphi}$) is a family of subsets of $\kappa$, $$\{E_t \colon t \in [\kappa]^{<\omega}\} \cup \{A_\alpha \colon \alpha < \kappa^+\}$$ satisfying the following conditions:
   \begin{itemize}
   	\item [(1)] $E_s\cap E_t = \emptyset$ for all $s, t \in [\kappa]^{<\omega}$ with $s \not =t$,
   	\item [(2)] $|\bigcap_{\alpha \in p} A_\alpha \cap  \bigcup_{t\not \supseteq \hat{\varphi}(p)}E_t| < \kappa$ for each $p \in [\kappa^+]^{<\omega}$,
   	\item [(3)] if $\hat{\varphi}(p) \subseteq t$, then $|\bigcap_{\alpha\in p} A_\alpha \cap E_t|=\kappa$ for each $p \in [\kappa^+]^{<\omega}$ and $t \in [\kappa]^{<\omega}$. 
   \end{itemize}
   
   Let $I$ be an index set and $\mathcal{F}$ be a filter on $\kappa$. The family $$\{E_t^i \colon t \in [\kappa]^{<\omega}, i \in I\} \cup \{A_\alpha^i \colon \alpha < \kappa^+, i \in i\}$$ is an \textit{independent matrix of $|I|$ step-families} (over $\kappa$) with respect to $\mathcal{F}, \hat{\varphi}$ iff
   \begin{itemize}
   	\item [(1)] for each fixed $i \in I$, $\{E_t^i \colon t \in [\kappa]^{<\omega}\} \cup \{A_\alpha^i \colon \alpha < \kappa^+\}$ is a step-family,
   	\item [(2)] if $n \in \omega, p_0, p_1, ..., p_{n-1} \in [\kappa^+]^{<\omega}, t_0, t_1, ..., t_{n-1} \in [\kappa]^{<\omega}$,  $i_0, i_1, ..., i_{n-1} \in I$ with $i_k\not = i_m, k\not = m$ and $\hat{\varphi}(p_k) \subseteq t_k$, then 
   	$$\bigcap_{k=1}^{n-1}(\bigcap_{\alpha\in p_k}A^{i_k}_{\alpha} \cap E^{i_k}_{t_k}) \in \mathcal{F}^+,$$
   	where $\mathcal{F}^+ = \{D\subseteq \kappa \colon \kappa\setminus D \not \in \mathcal{F}\}$.  
   \end{itemize}
   \noindent
   \textbf{Fact 1 (\cite{BK}).} If $\kappa$ is a regular cardinal and $\hat{\varphi}$ is a $\kappa$-shrinking function, then there exists and independent matrix of $2^\kappa$ step-families over $\kappa$ with respect to the  filter $\mathcal{RF}(\kappa)$, $\hat{\varphi}$.

   \section{Main result and open problem}
  
  The main result of this paper is Theorem 1 which immediately follows from Lemma 1 and Lemma 2.
  
  \begin{theorem}
  	Let $\kappa$ be a regular cardinal and let $\hat{\varphi}$ be a $\kappa$-shrinking function.
  	There exists a chain order isomorphic to $(2^\kappa)^+$ in the Rudin-Frol\'ik order of ultrafilters in $\beta \kappa$.
  \end{theorem}
  
   Lemma 1 can be proved with using an independent matrix of $2^\kappa$ step-families over $\kappa$ with respect to $\mathcal{RF}(\kappa), \hat{\varphi}$, but the proof of \cite[Theorem 2.1]{BE3}  can be adopted here as well. Since the first idea of the proof provides more complicated proof that the second one we decide to show the second one. This way we obtain more general lemma.
      
   \begin{lemma}
   	Let $\kappa, \tau$ be cardinals with $\kappa$ being regular. Let $$\{\mathcal{F}_{\xi, \zeta} \colon \xi< \tau, \zeta< \kappa, \mathcal{FR}(\kappa)\subseteq \mathcal{F}_{\xi, \zeta}\}$$ be a stratified set of filters. Then there exists a stratified set of ultrafilters $\{\mathcal{G}_{\xi, \zeta} \colon \xi< \tau, \zeta < \kappa\}$ such that $\mathcal{F}_{\xi, \zeta} \subseteq \mathcal{G}_{\xi, \zeta}$.	
   \end{lemma}
   
   \begin{proof}
   	We will proceed by induction on $\alpha$ in $2^\kappa$ steps.
   	Enumerate all subsets of $\kappa$ by
   	$$\{Z_\alpha \colon \alpha < 2^\kappa\}.$$
   	In step  $\alpha< 2^\kappa$  we wiil construct filters $\{\mathcal{F}^\alpha_{\xi, \zeta} \colon \xi < \tau, \zeta< \kappa\}$ fulfilling the following properties
   	\begin{itemize}
   		\item [(1)] $\mathcal{F}^0_{\xi, \zeta} = \mathcal{F}_{\xi, \zeta}$,
   		\item [(2)] $\{\mathcal{F}^\alpha_{\xi, \zeta} \colon \xi < \tau, \zeta< \kappa\}$ is a stratified set,
   		\item [(3)] $Z_\alpha$ or $\kappa\setminus Z_\alpha$ belongs to $\mathcal{F}^{\alpha+1}_{\xi, \zeta}$,
   		\item [(4)] $\mathcal{F}^{\alpha}_{\xi, \zeta} = \bigcup_{\eta< \alpha} \mathcal{F}^{\eta}_{\xi, \zeta}$ for $\alpha$ being a limit number.
   	\end{itemize}
   Assume that we have constructed the filters $\mathcal{F}^{\alpha}_{\xi, \zeta}$	 for each $\xi<\tau, \zeta< \kappa$.
   Now, we will show the successor step. Define
   $$W=\{\mathcal{F}^{\alpha}_{\xi, \zeta}, \colon Z_\alpha \in \mathcal{F}^{\alpha}_{\xi, \zeta}\}.$$
   Notice that if $\mathcal{F}^{\alpha}_{\xi, \zeta} \not \in \tilde{W}$ then $[\mathcal{F}^{\alpha}_{\xi, \zeta}, \{\kappa \setminus Z_\alpha\}]$ is not a filter. We will show that $[\mathcal{F}^{\alpha}_{\xi, \zeta}, \{Z_\alpha\}]$ is a filter, whenever $\mathcal{F}^{\alpha}_{\xi, \zeta} \in \tilde{W}$. 
   
   Indeed. Supose in contrario that 
   $$\mu = \min\{\gamma < \kappa^+ \colon \mathcal{F}^{\alpha}_{\delta, \varepsilon} \in W(\gamma)\setminus W \textrm{ and }[\mathcal{F}^{\alpha}_{\delta, \varepsilon}, \{Z_\alpha\}]\}$$
   is not a filter.
   Then by $(3)$, $\kappa \setminus Z_\alpha \in \mathcal{F}^{\alpha}_{\delta,\varepsilon}$.
   
   By $(2)$, the set
   $$\{\mathcal{F}^\alpha_{\psi, \theta} \colon \psi < \alpha, \theta < \kappa\}$$
   is a stratified set.
   Hence there exists $\rho > \xi$ such that the set
    $$\{\mathcal{F}^\alpha_{\rho, \lambda} \colon \mathcal{F}^\alpha_{\rho, \lambda} \in W(\eta) \textrm{ and } \kappa \setminus Z_\alpha \in \mathcal{F}^\alpha_{\rho, \lambda}\}$$
    is nonempty for $\eta < \gamma$.
    But $\gamma$ was assumed as a minimal one of this property. A contradiction.

    Thus, we can define
    $$\mathcal{F}^\alpha_{\xi, \zeta} = [\mathcal{F}^\alpha_{\xi, \zeta}, \{Z_\alpha\}]  \textrm{ whenever } \mathcal{F}^\alpha_{\xi, \zeta} \in \tilde{W}$$ and  $$\mathcal{F}^\alpha_{\xi, \zeta} = [\mathcal{F}^\alpha_{\xi, \zeta}, \{\kappa \setminus Z_\alpha\}] \textrm{ whenever } \mathcal{F}^\alpha_{\xi, \zeta} \not\in \tilde{W}$$
    
    To complete the proof is is enough to show that $\{\mathcal{F}^{\alpha+1}_{\xi, \zeta} \colon \xi<\alpha, \zeta < \kappa\}$ forms a stratified set.
    
    By $(2)$, and the definitions of $W$ and $\tilde{W}$ we have that 
    $$\{\mathcal{F}^{\alpha+1}_{\gamma, \delta} \colon A \cap Z_\alpha \in \mathcal{F}^{\alpha+1}_{\rho, \lambda}\}$$
    has cardinality $\kappa$, whenever $A \in \mathcal{F}^\alpha_{\xi, \zeta} \in \tilde{W}$ and $\rho>\xi$. 
    
    If $\mathcal{F}^\alpha_{\xi, \zeta}\not \in \tilde{W}$ then for each $\rho> \xi$ and for each $A \in \mathcal{F}^\alpha_{\xi, \zeta}$ the set
    $$\{\mathcal{F}^\alpha_{\rho, \lambda} \colon A \in \mathcal{F}^\alpha_{\rho, \lambda} \in \tilde{W}\}$$
    has cardinality $\kappa$. Hence
    $$\{\mathcal{F}^{\alpha+1}_{\rho, \lambda} \colon A \cap \{\kappa\setminus Z_\alpha\} \in \mathcal{F}^{\alpha+1}_{\rho,\lambda}\}$$
    also has cardinality $\kappa$.
    The successor step is complete.
    
    (It is easy to show that we can proceed our induction independently on $\mathcal{F}^\alpha_{\xi, \zeta}$ is a filter or an ultrafilter).
    
    Finally, put 
    $$\mathcal{G}_{\xi, \zeta} = \bigcup_{\alpha< 2^\kappa} \mathcal{F}^\alpha_{\xi, \zeta}.$$
    By $(2)$, the set
     $$\{\mathcal{G}_{\xi, \zeta} \colon \xi < \alpha, \zeta< \kappa\}$$ is stratified and by $(3)$, $\mathcal{G}_{\xi, \zeta}$ is an ultrafilter.
   	\end{proof}
   \\
   
   In the proof of Lemma 2 we use some ideas presented in the proof of \cite[Theorem 2.2]{BE3}.
   
   \begin{lemma}
   	Let $\kappa$ be a regular cardinal  and let $\hat{\varphi}$ be a $\kappa$-shrinking function. Then there exists a set of ultrafilters
   	$$\{\mathcal{G}^\alpha_{\xi, \zeta} \colon \xi < \alpha, \zeta < \kappa, \alpha < (2^\kappa)^+\}$$
   	fulfilling the following conditions
   	\begin{itemize}
   		\item $\{\mathcal{G}^\alpha_{\xi, \zeta} \colon \zeta < \kappa\}$ is a $\kappa$-discrete set of ultrafilters for each $\xi < \alpha < (2^\kappa)+$,
   		\item $\mathcal{G}^\alpha_{\xi, \zeta} <_{RF} \mathcal{G}^\beta_{\xi, \zeta}$ whenever  $\xi < \alpha < \beta < \tau$ and $\zeta < \kappa$.
   	\end{itemize}
   \end{lemma}
   
   \begin{proof} We will construct the proposed set of ultrafilters by the inductions on $\alpha$ in $(2^{\kappa})^+$ steps. At each step $\alpha< (2^\kappa)^+$ we will construct all ultrafilters $\mathcal{G}^\alpha_{\xi, \zeta}$ for all $\xi< \alpha$  and $\zeta< \kappa$. In order to do this we use a matrix of $2^\kappa$ step-families over $\kappa$ independent with respect to $\mathcal{FR}(\kappa), \hat{\varphi}$.
   	
   		By Fact 1, fix a matrix $$\{E_t^i \colon t \in [\kappa]^{<\omega}, i \in 2^\kappa\} \cup \{A_\eta^i \colon \eta < \kappa^+, i \in 2^\kappa\}$$ of $2^\kappa$ step-families over $\kappa$ independent with respect to the filter $\mathcal{FR}(\kappa), \hat{\varphi}$. 
   		Assume that for each $i \in 2^\kappa$ we have 
   		\begin{itemize}
   			\item [(a)] $E^i_s\cap E^i_t = \emptyset$ for all $s, t \in [\kappa]^{<\omega}$ with $s \not =t$,
   			\item [(b)] $|\bigcap_{\eta \in p} A^i_\eta \cap  \bigcup_{t\not \supseteq \hat{\varphi}(p)}E^i_t| < \kappa$ for each $p \in [\kappa^+]^{<\omega}$,
   			\item [(c)] if $\hat{\varphi}(p) \subseteq t$, then $|\bigcap_{\eta\in p} A^i_\eta \cap E^i_t|=\kappa$ for each $p \in [\kappa^+]^{<\omega}$ and $t \in [\kappa]^{<\omega},$ 
   			\item [(d)] $\bigcup\{E^i_t \colon t \in [\kappa]^{< \omega}\} = \kappa,$
   			\item [(e)] $|\bigcap_{\eta \in p} A^i_\eta \setminus  \bigcup_{t \supseteq \hat{\varphi}(p)}E^i_t| < \kappa$ for each $p \in [\kappa^+]^{<\omega}.$
   		\end{itemize} 
   	
   	Note that $(b)$ and $(d)$ implies $(e)$.
   	Moreover, the condition $(b)$ is still preserved after expanding $\{E_t^i \colon t \in [\kappa]^{<\omega}\}$ to a partition of $\kappa$. 
   	\\Indeed. If there are $p_0 \in [\kappa^+]^{<\omega}$ and $i_0 \in 2^\kappa$ such that
   	$$|\bigcap_{\eta \in p_0} A_\eta^{i_0} \cap  \bigcup_{t\not \supseteq \hat{\varphi}(p_0)}E_t^{i_0}| = \kappa,$$
   	then $|\bigcup_{t\not \supseteq \hat{\varphi}(p_0)}E_t^{i_0}| = \kappa.$
   	Then, by $(d)$ and $(a)$, there would exist  $t_0 \supseteq \hat{\varphi}(p_0)$ such that  $|E_{t_0}^{i_0}| < \kappa$. Hence
   	$$|\bigcap_{\eta\in p_0} A^{i_0}_\eta \cap E^{i_0}_{t_0}|<\kappa.$$
   	which contradicts $(c)$. 
   	
   		 Now, we will define two families of sets
   		$$\{\{W^\alpha_{\gamma, \zeta}\colon \zeta< \kappa\} \colon \gamma < \alpha < (2^\kappa)^+\}$$
   		and
   		$$\{\{Z^\alpha_{\gamma, \zeta}\colon \zeta< \kappa\} \colon \gamma < \alpha < (2^\kappa)^+\}.$$
   		For each $\alpha < (2^\kappa)^+$ take families 
   		$$\{\{A^i_\eta \colon \eta < \kappa^+\} \colon i \in |\alpha|\}$$
   		and partitions
   		$$\{\{E^i_t \colon t \in [\kappa]^{< \omega}\} \colon i \in |\alpha|\}$$
   		and renumerate them by taking a bijection
   		$$r \colon \alpha \to |\alpha|$$
   		and increasing functions
   		$$g \colon (2^\kappa)^+ \to \kappa^+$$
   		and 
   		$$h \colon (2^\kappa)^+ \to \{t \colon t \in [\kappa]^{< \omega}\}, $$
   		i.e. the functions of the following properties: for all $\alpha< \beta < (2^\kappa)^+$  we have
   		$g(\alpha)<g(\beta)$ and if $h(\alpha) = t$ and $h(\beta) = s$ then $t \subset s$.
   		
   		Then we take
   		$$W^\alpha_{\gamma, \zeta} = A^{r(\alpha)}_{g(\gamma), \zeta}
   		\textrm{ and }
   		Z^\alpha_{\gamma, \zeta} = E^{r(\alpha)}_{h(\gamma), \zeta}.$$
   		Now, we will construct a family
   		$$\{B^\alpha_{\xi, \zeta}\colon \zeta< \kappa\}$$
   		of partitions of $\kappa$ such that $|B^\alpha_{\xi, \zeta}| = \kappa$ for each $\xi< \alpha< (2^\kappa)^+$, which will be used to the construction of required ultrafilters.
   		For this purpose we consider two functions
   		$$\Psi \colon (2^\kappa)^+ \to \{p \colon p \in [\kappa^+]^{<\omega}\}$$
   		$$\Phi \colon (2^\kappa)^+ \to \{t \colon t \in [\kappa]^{<\omega}\}.$$
   		and define $B^\alpha_{\xi, \zeta}$ as follows
   		$$B^\alpha_{\xi, \zeta} = \left\{\begin{array}{rcl}
   		(\bigcap_{\eta \in \Psi(\xi)}W^\alpha_{\eta, \zeta}) \setminus (\bigcup_{\sigma \in \Psi(\xi)}Z^\alpha_{\sigma, \zeta}) & \hat{\phi} (\Psi(\xi)) \subseteq \Phi(\xi)\\
 (\bigcap_{\eta \in \Psi(\xi)}W^\alpha_{\eta, \zeta}) \cap (\bigcup_{\sigma \in \Phi(\xi)} Z^\alpha_{\sigma, \zeta}) & \hat{\phi} (\Psi(\xi)) \not\subseteq \Phi(\xi)
   		\end{array} \right.$$
   		Now, we are ready to start the construction of the required ultrafilters.
   		
   		Let $\{\mathcal{G}^1_{0, \zeta} \colon \zeta < \kappa\}$ be a family of arbitrary sets of ultrafilters such that
   		$$[\mathcal{FR}(\kappa), \{B^1_{0, \zeta}\}] \subseteq \mathcal{G}^1_{0, \zeta}$$
   		for any $\zeta< \kappa$.
   		  
  In order to define $\mathcal{G}^\beta_{\delta, \zeta}$, where $\beta = \alpha + 1$, we define a filter
  $$\mathcal{F}_{\beta, \eta} = [\mathcal{FR}(\kappa), \{B^\alpha_{\xi, \zeta} \colon \xi< \gamma, \lambda \in B^\gamma_{\xi, \zeta}\}, \{B^\alpha_{\gamma, \lambda}\}, \bigcup_{\gamma< \xi< \alpha}\{\bigcup_{\xi \in D} B^\alpha_{\xi, \zeta} \colon D \in \mathcal{G}^\xi_{\gamma, \lambda}\}]$$
  for each $\gamma< \alpha$ and $\lambda< \kappa$. 
  
  Notice that $\mathcal{F}_{\gamma, \lambda}$ are filters because we only use the elements of the matrix odf $2^\kappa$ step-families over $\kappa$ independent with respect to $\mathcal{FR}(\kappa), \hat{\varphi}$ for their construction.
  
  It follows from the construction that $\{\mathcal{F}_{\gamma, \lambda} \colon \lambda < \kappa\}$ is a $\kappa$-discrete set for each $\gamma < \alpha$.
 Now, we will show that the second condition in Lemma 2 is also fulfilled.
Let $A \in \mathcal{F}_{\gamma, \lambda}$ and $\delta > \gamma$.
 
Analyzing the construction of $\mathcal{F}_{\gamma, \lambda}$ it is enough to show that 
$$A \supseteq B^\alpha_{\xi_1, \eta}\cap B^\alpha_{\gamma, \lambda} \cap \bigcap_{i=2}^{4} \bigcup_{j \in D_i} B^\alpha_{\xi_i, j},$$
where  $\xi_1 < \gamma < \xi_2<\xi_3 = \delta < \xi_4$ and $D_i \in \mathcal{G}^{\xi_i}_{\gamma, \lambda}.$ (The proof of the other forms of elements of $\mathcal{F}_{\gamma, \lambda}$ runs in the similar way but it is only more complicated to describe).
By the above construction we have
\begin{itemize}
\item [(1)] $B^\alpha_{\xi_1, \eta} \in \mathcal{F}_{\delta, \lambda}$, whenever
$\lambda \in B^\delta_{\xi_1, \eta} \in \mathcal{G}^\delta_{\xi_1, \eta}$.
\item [(2)] $B^\delta_{\xi_1,\eta} \in \mathcal{G}^\delta_{\gamma, \nu}$ such that $\nu \in B^\gamma_{\xi_1, \eta}.$
\item [(3)]  $B^\delta_{\xi_1,\eta} \in \mathcal{G}^\delta_{\gamma, \lambda}$.
\item [(4)] $B^\alpha_{\gamma, \lambda} \in \mathcal{F}_{\delta, \zeta}$, whenever
$\zeta \in B^\delta_{\gamma, \lambda} \in \mathcal{G}^\delta_{\gamma, \lambda}$.
\item [(5)] $\bigcup_{j \in D_2} B^\alpha_{\xi_2, j}\in \mathcal{F}_{\delta,\zeta},$
whenever $D_2 \in \mathcal{G}^{\xi_2}_{\gamma,\lambda}$ and $\zeta \in \bigcup_{j \in D_2}B^\delta_{\xi_2, j}$.
\item [(6)] $\bigcup_{j \in D_3}B^\alpha_{\delta,j} \in \mathcal{F}_{\delta, \zeta}$, whenever $\zeta \in D_3 \in \mathcal{G}^\delta_{\gamma, \lambda}$.
\item [(7)] $D_4 \in \mathcal{G}^{\xi_4}_{\gamma, \lambda}$ which is the consequence of $\mathcal{G}^\delta_{\gamma, \lambda} <_{RK} \mathcal{G}^{\xi_4}_{\gamma, \lambda}$. Hence there exist $C \in \mathcal{G}^\delta_{\gamma, \lambda}$ and $G_\eta \in \mathcal{G}^{\xi_4}_{\delta, \eta}$ such that $D_4 = \bigcup_{\eta\in C} Z_\eta$. Then $\bigcup_{j \in Z_\eta} B^\alpha_{\xi_4, j} \in \mathcal{F}_{\delta, \eta}$ for each $\eta \in C$ and $\bigcup_{j \in D_4} B^\alpha_{\xi_4, j} \in \mathcal{F}_{\delta, \eta}$. Thus $\bigcup_{j \in D_4} B^\alpha_{\xi_4, j} \in \mathcal{F}_{\delta, \zeta}$ for each $\zeta \in C$.
\end{itemize} 

From $(1)-(7)$ follows that 
$A \in \mathcal{F}_{\delta, \zeta}$ whenever $\zeta \in B^\delta_{\xi_1, \eta} \cap B^\delta_{\gamma, \lambda} \cap \bigcup_{j \in D_2}B^\delta_{\xi_2, j} \cap D_3\cap C$. Hence $A \in \mathcal{G}^\delta_{\gamma, \lambda}$.

By Lemma 1, there exists a stratified set of ultrafilters 
 $$\{\mathcal{G}^\alpha_{\gamma, \lambda} \colon \gamma < \alpha, \lambda<\kappa\}$$
 such that $\mathcal{G}^\alpha_{\gamma, \lambda} \supseteq \mathcal{F}_{\gamma, \lambda}$. By the above construction
  $$\mathcal{G}^\delta_{\gamma, \lambda} \subseteq \{\{\eta \colon A \in \mathcal{F}_{\delta, \eta}\} \colon A \in \mathcal{F}_\gamma, \lambda\}.$$ Hence $\mathcal{G}^\delta_{\gamma, \lambda} <_{RF} \mathcal{G}^\alpha_{\gamma, \lambda}$.
   		\end{proof}

   \noindent
   \textbf{Open problem} Is there a chain order isomorphic co $(2^\kappa)^+$ in the Rudin-Frol\'ik order of ultrafilters in $\beta\kappa$ for $\kappa$-singular?
   \\
  \\ 
   \textbf{Conflict of interest} There is no conflict of interest.
   \\
   \\
   \textbf{Data availability} Not applicable.

		{\sc Joanna Jureczko}
		\\
		Wroc\l{}aw University of Science and Technology, Wroc\l{}aw, Poland
		\\
		{\sl e-mail: joanna.jureczko@pwr.edu.pl}


\begin {thebibliography}{123456}
	\thispagestyle{empty}
	\bibitem{BK} J. Baker, K. Kunen, Limits in the uniform ultrafilters. Trans. Amer. Math. Soc. 353 (2001), no. 10, 4083–4093.
	
	\bibitem{DB} D. Booth, Ultrafilters on a countable set, Ann. Math. Logic 2 (1970/71), no. 1, 1--24.
	
		\bibitem{BB} L. Bukovsk\'y, E. Butkovi\v cov\'a, Ultrafilters with $\aleph_0$ predecessors in Rudin-Frol\'ik order, Comment. Math. Univ. Carolin. 22 (1981), no. 3, 429–-447.
		
	\bibitem{BE1} E. Butkovi\v cov\'a, Ultrafilters without immediate predecessors in Rudin-Frolík order. Comment. Math. Univ. Carolin. 23 (1982), no. 4, 757–-766.
	
	\bibitem{BE3} E. Butkovi\v cov\'a, Long chains in Rudin-Frol\'ik order, Comment. Math. Univ. Carolin. 24 (1983), no. 3, 563–-570.
	
	\bibitem{BE4} E. Butkovi\v cov\'a, Subsets of $\beta\mathbb{N}$ without an infimum in Rudin-Frol\'ik order, Proc. of the 11th Winter School on Abstract Analysis, (Zelezna Ruda 1983), Rend. Circ. Mat. Palermo (2) (1984), Suppl. no. 3, 75--80.
	
	\bibitem{BE5} E. Butkovi\v cov\'a, Decrasing chains without lower bounds in the Rudin-Frol\'ik order, Proc. AMS, 109, (1990) no. 1, 251--259.
		\bibitem{CN} W. W. Comfort, S. Negrepontis, The Theory of Ultrafilters, Springer 1974.
		
		\bibitem{ZF} Z. Frol\'ik, Sums of ultrafilters. Bull. Amer. Math. Soc. 73 (1967), 87--91.
		
			\bibitem{MG} M. Gitik, Some constructions of ultrafilters over a measurable cardinal, Ann. Pure Appl. Logic 171 (2020) no. 8, 102821, 20pp.
			
				\bibitem{TJ}    Jech, T., Set Theory, The third millennium edition, revised and expanded. Springer Monographs in Mathematics. Springer-Verlag, Berlin, 2003.
		
	\bibitem{JJ_order} J. Jureczko, Ultrafilters without immediate predecessors in Rudin-Frol\'ik order for regulars, Results Math. 77, 230 (2022).
	
	\bibitem{JJ_order4} J.Jureczko, Chians in the Rudin-Frol\.ik order for regulars, Raporty Katedry Telekomunikacji i Teleinformatyki. 2022, Ser. PRE nr 25, 12 s., (access: dona.pwr.edu.pl)
			
	\bibitem{JJ_order2} J. Jureczko, A note on special subsets of the Rudin-Frol´ık order	for regulars,  https://arxiv.org/pdf/2304.02143.pdf, (accepted in Math. Slovaca 19.09.2022).
	
	\bibitem{JJ_order3} J. Jureczko, Decreasing chains without lower bounds in the Rudin Frol\.ik order for regulars, https://arxiv.org/pdf/2304.01398.pdf. 
	
		\bibitem{JJ_kanamori} J. Jureczko, On some constructions of ultrafilters over a measurable cardinal, (in preparation).
		
			\bibitem{AK} A. Kanamori, Ultrafilters over a measurable cardinal, Ann. Math. Logic, 11 (1976), 315--356.
	
	\bibitem{KK} K. Kunen, Weak P-points in $\mathbb{N}^*$. Topology, Vol. II (Proc. Fourth Colloq., Budapest, 1978), pp. 741–749, Colloq. Math. Soc. János Bolyai, 23, North-Holland, Amsterdam-New York, 1980.
	
		\bibitem{MR1} M.E. Rudin, Types of ultrafilters in: Topology Seminar Wisconsin, 1965 (Princeton Universiy Press, Princeton 1966).
	
	\bibitem{MR} M. E. Rudin, Partial orders on the types in $\beta \mathbb{N}$. Trans. Amer. Math. Soc. 155 (1971), 353--362.
	
		\end{thebibliography}
\end{document}